\documentclass[12pt,pagebackref]{article}
\usepackage{amsfonts,amsmath,amsthm,yfonts,mathdots,
latexsym,
amssymb 
}
\newtheorem{theorem}{Theorem}
\newtheorem{lemma}[theorem]{Lemma}

\def\ex{{\rm ex}}

\title{Extremal numbers for odd cycles}
\author{Zoltan F\"uredi${}^\star$ \enskip and \enskip David S. Gunderson${}^\dagger$}
\date{${}^\star$ R\'enyi Institute of Mathematics, 
  Hungarian Academy of Sciences
  \\
E-mail: {\tt z-furedi@illinois.edu, furedi.zoltan@renyi.mta.hu}
\footnotetext{\noindent${}^\star$
Research supported in part by the Hungarian National Science Foundation OTKA 104343,
and by the European Research Council Advanced Investigators Grant 267195.\\
\indent\indent${}^\dagger$
 Research supported by NSERC Discovery grant  228064. \\
MSC-class: 05C35, 05D99              \hfill  \jobname\quad  24 October 2013
 \hfill Printed on {\today}
}\\
 \medskip ${}^\dagger$ University of Manitoba, Winnipeg, Canada\\
E-mail: {\tt David.Gunderson@umanitoba.ca}
}

\begin{document}
\maketitle

\begin{abstract}
We describe the $C_{2k+1}$-free graphs on $n$ vertices with maximum number of edges.
The extremal graphs are unique for $n\notin \{ 3k-1, 3k, 4k-2, 4k-1\}$.
The value of $\ex(n,C_{2k+1})$ can be read out from the works of Bondy~\cite{Bondy2},
 Woodall~\cite{Woodall:72}, and Bollob\'as~\cite{BB:79}, but here we give a new
 streamlined proof.
The complete determination of the extremal graphs is also new.

We obtain that the bound for $n_0(C_{2k+1})$ is $4k$ in the classical theorem of Simonovits,
 from which the unique extremal graph is the bipartite Tur\'an graph.
\end{abstract}

\newpage
\section{Introduction, exact Tur\'an numbers}
Given a class of simple graphs ${\mathcal F}$ let us call a graph ${\mathcal F}$-\emph{free}
 if it contains no copy of $F$ as a (not necessarily induced) subgraph for each $F\in {\mathcal F}$.
Let $\ex(n; {\mathcal F })$
 denote the maximal number of edges in an ${\mathcal F}$-free graph on $n$ vertices.
If the class of graphs ${\mathcal F}=\{F_1, F_2, \dots \}$ consists of a single graph then
 we write $\ex(n;F)$ instead of $\ex(n;\{ F\})$.

Let $T_{n,p}$ denote the \emph{Tur\'an graph},  the complete equi-partite graph,
 $K_{n_1, n_2, \dots, n_p}$ where $\sum _i n_i=n$ and
 $\lfloor n/p \rfloor\leq  n_i\leq \lceil n/p \rceil$. 
By Tur\'an's theorem~\cite{Turan:41, Turan:54} we have $\ex(n; K_{p+1})=e(T_{n,p})$;
 furthermore, $T_{n,p}$
 is the unique $K_{p+1}$-free graph that attains the extremal number.
The case  $\ex(n; K_3)=\lfloor n^2/4\rfloor$ was shown earlier by
 Mantel~\cite{Mantel:07}.

There are very few cases when the Tur\'an number
  $\ex(n;{\mathcal F})$ is known exactly for all $n$.
One can mention the case when $F=M_{\nu+1}$ is a matching of a given size, $\nu+1$.
Erd\H os and Gallai~\cite{EG:59} showed that
$$
   \ex(n, M_{\nu+1})= \max\{ {2\nu+1 \choose 2}, {\nu\choose 2}+\nu(n-\nu)\}.
   $$
For the path of $k$ vertices Erd\H os and Gallai~\cite{EG:59}
 proved an asymptotic and $\ex(n;P_k)$ was determined  for all $n$ and $k$ by
Faudree and Schelp~\cite{FaudreeSchelpJCT75} and independently by Kopylov~\cite{Kopylov}.
Erd\H os and Gallai~\cite{EG:59} proved an asymptotic
 for the class of long cycles ${\mathcal C}_{\ge \ell}:= \{ C_\ell, C_{\ell+1}, C_{\ell+2}, \dots\}$.
The exact value of the Tur\'an number  $\ex(n;{\mathcal C}_{\ge \ell})$
  was determined by Woodall~\cite{WoodallA}
  and independently and at the same time by Kopylov~\cite{Kopylov}.

There is one outstanding result which gives infinitely many exact Tur\'an numbers,
 Simonovits' chromatic critical edge theorem~\cite{SimSymm77}.
It states that if $\min \{ \chi(F): F\in {\mathcal F} \}=p+1\geq 3$ and
 there exists an $F\in {\mathcal F}$ with an edge $e\in E(F)$ such that by removing
 this edge one has $\chi(F-e)\leq p$, then there exists an $n_0({\mathcal F})$
 such that  $T_{n,p}$ is the only extremal graph for ${\mathcal F}$ for $n \geq n_0$.
The authors are not aware of any (non-trivial) further result when $\ex(n, {\mathcal F})$
 is known for all $n$, neither any $F$ for which the value of $n_0(F)$ had been determined,
  except the case of odd cycle discussed below.

\section{The result, the extremal graphs without  $C_{2k+1}$}

The aim of this paper is to determine the Tur\'an number of odd cycles for all $n$ and $C_{2k+1}$
 together with the extremal graphs.
The value of $\ex(n,C_{2k+1})$ can be read out from the works of Bondy~\cite{Bondy:71,Bondy2},
 Woodall~\cite{Woodall:72}, and Bollob\'as~\cite{BB:79} (pp. 147--156)
 concerning (weakly) pancyclic graphs.
For a recent presentation see Dzido~\cite{Dzido} who also considered the Tur\'an number of wheels.
But here we give a new streamlined proof and a complete description of the extremal graphs.

Since  $K_{\lceil n/2\rceil, \lfloor n/2 \rfloor}$ contains no odd cycles,
for any $k\geq 1$, $\ex(n; C_{2k+1})\geq \lfloor n^2/4\rfloor$.
For $C_3$ here equality holds for all $n$ with the only extremal graph is $T_{n,2}$
 by the Tur\'an-Mantel's theorem.
From now on, we suppose that $2k+1\geq 5$.
Also for $n\leq 2k$ obviously $\ex(n, C_{2k+1})={n \choose 2}$ so we may suppose
 that $n\geq 2k+1$.

Every edge of an odd cycle is color critical so Simonovits' theorem implies that
 the complete bipartite graph is the only extremal graph and
 $\ex(n; C_{2k+1})= e(T_{n,2})= \lfloor n^2/4\rfloor$  for $n \geq n_0(C_{2k+1})$.
After choosing the right tools we present a streamlined proof and show
 that $n_0(C_{2k+1})=4k$ (in case of $2k+1\geq 5$).

We define two classes of $C_{2k+1}$-free graphs which could have at least as
 many edges as $T_{n,2}$ for $n\leq 4k-1$.
A {\it cactus} $B(n; n_1, \dots, n_s)$ (for $n\geq 2$, $s\geq 1$ with $\sum_i (n_i-1)=n-1$)
 is a connected graph where the $2$-connected blocks are complete graphs of sizes $n_1, \dots, n_s$.
Let us denote by $g(n,k)$ the largest size of an $n$-vertex cactus avoiding $C_{2k+1}$.
For this maximum all block sizes should be exactly $2k$ but at most one which is smaller.
Write $n$ in the form $n=(s-1)(2k-1)+r$ where $s\geq 1$, $2\le r\le 2k$ are integers. Then
\begin{equation}\label{eq:gnk}
 g(n,k)=(s-1){2k\choose 2} + {r\choose 2}. \end{equation}
Note that $g(n,k)> \lfloor n^2/4\rfloor$ for $3\leq n\leq 4k-3$ and we have
 $g(n,k)=e(T_{n,2})= \lfloor n^2/4\rfloor$ if $n\in \{ 4k-2, 4k-1\}$.
Thus  the Simonovits threshold  $n_0(C_{2k+1})$ is at least $4k$.

For $n\geq k$, define the graph $H_1(n,k)$ on $n$ vertices by its degree sequence; it has $k$
vertices of degree $n-1$ and all other vertices have degree $k$. 
Then $H_1(n,k)$ is a complete bipartite graph $K_{k, n-k}$, together with all possible edges added
in the first partite set.  This graph does not contain the cycle $C_{2k+1}$.
Letting $h_1(n,k)$ denote the size of $H_1(n,k)$,
\begin{equation}\label{eq:h1nk}
  h_1(n,k)={k\choose 2} + k(n-k). \end{equation}
Note that $h_1(n,k)\leq g(n,k)$ for all $k\leq n$ and here
 equality holds if  $n$ is in the form $n=(s-1)(2k-1)+r$ where $s\geq 1$ and
 $r\in \{k, k+1\}$.

 \begin{theorem}\label{th:main}
 For any  $n\geq 1$ and $2k+1\geq 5$,
 $$
\ex(n;C_{2k+1})=
\begin{cases} {n \choose 2} & \text{ for }\,  n \leq 2k,  \\
g(n,k) & \text{ for }\, 2k+1\leq n\leq 4k-1 \, \text{ and}\\
 \lfloor n^2/4\rfloor & \text{ for }\,  n\geq 4k-2.
\end{cases}$$
Furthermore, the only extremal graphs are $K_n$ for $n\leq 2k$; $B(n; 2k, n-2k+1)$
 for $2k+1\leq n\leq 4k-1$; $H_1(n,k)$ for $n \in \{ 3k-1, 3k\}$; and
 the complete bipartite graph $K_{\lceil n/2\rceil, \lfloor n/2 \rfloor}$ for $n\geq 4k-2$.
  \end{theorem}

 \section{A lemma on $2$-connected graphs without $C_{2k+1}$}\label{S:lemma}

\begin{lemma}\label{le:main}
 Suppose that $n\geq 2k+1\geq 5$ and $G$ is a $2$-connected, $C_{2k+1}$-free, non-bipartite graph
with at least $\lfloor n^2/4\rfloor$ edges.
Then $e(G)\leq \ex(n; C_{2k+1})$ and here equality holds
 only if $n\in  \{ 3k-1, 3k\}$ and $G=H_1(n,k)$.
  \end{lemma}

For $5\leq 2k+1\leq n$, define the graph $H_2(n,k)$ on $n$ vertices
 and
$$ h_2(n,k):= {2k-1\choose 2}+ 2(n-2k+1)
 $$
edges,  consisting of a complete graph $K_{2k-1}$ containing two special
vertices which are connected to all other vertices. Then $H_2(n,k)$ is a
$2$-connected $C_{2k+1}$-free graph.
For $k=2$ the graphs $H_1(n,k)$ and $H_2(n,k)$ are isomorphic.
Recall a result of Kopylov~\cite{Kopylov} in a form we use it:
Suppose that the $2$-connected graph $G$ on $n$ vertices
 contains no cycles of length $2k+1$ or larger and $n\geq 2k+1\geq 5$.
Then
\begin{equation}\label{eq:Kop}
 e(G)\leq \max\{ h_1(n,k), h_2(n, k)\}
  \end{equation}
and this bound is the best possible.
Moreover, only the graphs $H_1(n,k)$ and $H_2(n,k)$ could be extremal.
For further explanation and background see the recent survey~\cite{FureSimoErdos100}.

The other result we need is due to Brandt~\cite{Brandt1997}:
Let $G$ be a non-bipartite graph of order $n$  and suppose that
\begin{equation}\label{eq:Brandt}
 e(G)> (n-1)^2/4+1,
  \end{equation}
then $G$ contains cycles of every length between $3$ and the length of its longest cycle.

\medskip
\noindent{\bf Proof} of Lemma~\ref{le:main}:\enskip
The inequality $e(G)\leq \ex(n, C_{2k+1})$ follows from
 the definition.
Suppose that here equality holds.
Apply Brandt's theorem~\eqref{eq:Brandt}.
We obtain that $G$ contains cycles of all lengths $3,4, \dots, \ell$ where $\ell$ stands for the
 longest cycle length in $G$.
It follows that $\ell\leq 2k$.
Kopylov's theorem~\eqref{eq:Kop} implies that
\begin{equation*}
   \max\{ g(n,k), \lfloor n^2/4\rfloor\}\leq \ex(n, C_{2k+1})=e(G)\leq \max\{ h_1(n,k), h_2(n, k)\}.
  \end{equation*}
Since $g(n,k)> h_2(n,k)$ except for $(n,k)\in \{ (5,2), (6,2)\}$ and
$g(n,k)> h_1(n,k)$ except if  $n$ is in the form $n=(s-1)(2k-1)+r$ where $s\geq 2$ and
 $r\in \{k, k+1\}$ we obtain that $e(G)=h_1(n,k)$,
 $n$ should be in this form, and $G=H_1(n,k)$.

Finally, $h_1(n,k)< \lfloor n^2/4\rfloor$ for $n\geq 4k$ so we obtain that
 indeed $n\in \{ 3k-1, 3k\}$.
\qed

\section{The proof of Theorem~\ref{th:main}}

Suppose that $G$ is an extremal $C_{2k+1}$-free graph, $e(G)=\ex(n, C_{2k+1})$.
Then $G$ is connected.
Consider the cactus-like block-decomposition of $G$, $V(G)=V_1\cup V_2\cup \dots \cup V_s$,
 where the induced subgraphs $G[V_i]$ are either edges or maximal $2$-connected subgraphs of $G$.
Let $n_i:=|V_i|$, we have $n-1=\sum_i(n_i-1)$, and each $n_i\geq 2$.
We have $e(G[V_i])= \ex(n_i, C_{2k+1})$ otherwise one can replace $G[V_i]$ by an extremal graph
 of the same order $n_i$ and obtain another $C_{2k+1}$-free graph of size larger than $e(G)$.
Therefore $e(G[V_i])\geq \lfloor n_i^2/4\rfloor$ and there are three types of blocks\\
--- complete graphs (if $n_i\leq 2k$),\\
--- bipartite blocks
   with $e(G[V_i])=\lfloor n_i^2/4\rfloor$.
Finally,\\
--- if $n_i\geq 2k+1$ and $G[V_i]$ is not bipartite then Lemma~\ref{le:main} implies that
  $n_i\in \{ 3k-1, 3k\}$ and $G[V_i]=H_1(n_i,k)$.

We may rearrange the graphs $G[V_i]$ and the sets $V_i$ such a way that they share a common vertex
 $v\in \cap V_i$ and otherwise the sets $V_i\setminus \{ v\}$ are pairwise disjoint.
The obtained new graph $G^*$ also $C_{2k+1}$-free and extremal, it has the same size and order as $G$ has.

If $s=1$ then we are done.
Suppose $s\geq 2$.
If all blocks are complete graphs, then $e(G)\leq g(n,k)$.
Since $g(n,k)< e(T_{n,2})$ for $n> 4k-1$ we get that $n\leq 4k-1$
 and $G^*$ (and $G$) has only two blocks and at least one of them is of size $2k$.

Finally, suppose that there are two blocks $V_i$ and $V_j$, $|V_i|=a$ and $|V_j|=b$,
 such that $G[V_i]$ and $G[V_j]$ are not both complete subgraphs.
We claim that in this case one can remove
 the edges of $G[V_i]$ and $G[V_j]$ from $G^*$ and place
 a copy of $T_{a+b-1,2}$ or some other graph $L$ onto $V_i\cup V_j$
 such that the obtained new graph is $C_{2k+1}$-free and
 it has more edges than $e(G)$, a contradiction.

Indeed, if $G[V_i]$ is a large bipartite graph,
 $a:=n_i\geq 2k+1$, $G[V_i]=T_{a,2}$
 and  $G[V_j]$  is a complete bipartite graph, too,
 then we can increase $e(G^*)$ since
 \begin{equation}\label{eq:51}
  e(T_{a,2})+ e(T_{b,2})\leq \frac{1}{4}a^2+ \frac{1}{4}b^2<
   \lfloor \frac{1}{4}(a+b-1)^2\rfloor = e(T_{a+b-1,2}).
  \end{equation}
In the remaining cases the inequalities concerning the number of edges of $e(L)$
 are just elementary high school algebra.
If $G[V_i]=T_{a,2}$ and $G[V_j]=H_1(b,k)$ or $K_b$ then we can replace them again
 by a complete bipartite graph $T_{a+b-1,2}$.
From now on, we may suppose that each block is either a complete graph (of size at most $2k$)
 or an $H_1(a,k)$.
If $G[V_i]=H_1(a,k)$ for some $a\in \{ 3k-1, 3k\}$  and  $G[V_j]=H_1(b,k)$ (with $b\in \{ 3k-1, 3k\}$)
 or  $G[V_j]=K_b$ with $k\leq b\leq 2k$ then we replace $G[V_i]\cup G[V_j]$ again by
 a $T_{a+b-1,2}$.
Finally, if $G[V_i]=H_1(a,k)$ for some $a\in \{ 3k-1, 3k\}$  and  $G[V_j]=K_b$ with
 $2\leq b\leq k$ then we replace $G[V_i]\cup G[V_j]$ by two complete graphs of sizes
  $2k$ and $a+b-2k$
  and use  $e(H_1(a,k))+ e(K_b) <  e(B(a+b-1; 2k, a+b-2k))$
   to get a contradiction.
This completes the proof of the claim and Theorem~\ref{th:main}. \qed



\begin{thebibliography}{99} 



\parskip =2 pt plus 2pt
\small{

\bibitem{BB:79} B.~Bollob\'as,
{\it Extremal Graph Theory},
 Academic Press, New York, 1979.

\bibitem{Bondy:71} J. A. Bondy,
Pancyclic graphs I,
{\em J. Combin. Theory Ser. B} {\bf 11} (1971), 80--84.

\bibitem{Bondy2} J. A. Bondy,
Large cycles in graphs, {\em Discrete Math.} {\bf 1} (1971), 121--132.

\bibitem{Brandt1997}
S. Brandt, 
A sufficient condition for all short cycles,
{\it in:} 4th Twente Workshop on Graphs and Combinatorial Optimization (Enschede, 1995).
{\em Discrete Appl. Math.} {\bf 79} (1997), 63--66.

\bibitem{Dzido}
T. Dzido, A note on Tur\'an numbers for even wheels,
{\em  Graphs Combin.} {\bf 29} (2013), 1305--1309.

\bibitem{EG:59}P. Erd\H{o}s and T. Gallai,
On maximal paths and circuits of graphs,
{\em Acta Math. Acad. Sci. Hung.} {\bf 10} (1959), 337--356.

\bibitem{FaudreeSchelpJCT75}
R. J. Faudree and R. H. Schelp,
Path Ramsey numbers in multicolorings,
{\em J. Combin. Theory Ser. B} {\bf 19} (1975), 150--160.

\bibitem{FureSimoErdos100}
Z. F\"uredi and  M. Simonovits,
The history of degenerate (bipartite) extremal graph problems,
\emph{Bolyai Soc. Studies} (The Erd\H os Centennial) {\bf 25} (2013), 167--262.

\bibitem{Kopylov}
G. N. Kopylov,
Maximal paths and cycles in a graph,
{\em Dokl. Akad. Nauk SSSR} {\bf 234} (1977), no. 1, 19--21.
(English translation: {\em Soviet Math. Dokl.} {\bf 18} (1977), no. 3, 593--596.)

\bibitem{Mantel:07}
W. Mantel, 
Solution to Problem 28, by H. Gouwentak, 
 W. Mantel, J. Teixeira de Mattes, 
F. Schuh, 
and W. A. Wythoff, 
{\em Wiskundige Opgaven} {\bf 10} (1907), 60--61.

\bibitem{SimSymm77}
M. Simonovits,
Extremal graph problems with symmetrical extremal graphs, additional chromatic conditions,
{\em Discrete Math.} {\bf 7} (1974), 349--376.

\bibitem{Turan:41}
P. Tur\'an,
Eine Extremalaufgave aus der Graphentheorie (in Hungarian),
{\em Math. Fiz Lapook} {\bf 48} (1941), 436--452.

\bibitem{Turan:54}
P. Tur\'an,
On the theory of graphs, {\em Colloq. Math.} {\bf 3} (1954), 19--30.

\bibitem{Woodall:72}
D. R. Woodall,
Sufficient conditions for circuits in graphs,
{\em Proc. London Math. Soc.} (3) {\bf 24} (1972), 739--755.

\bibitem{WoodallA}
D. R. Woodall,
Maximal circuits of graphs I,
{\em Acta Math. Acad. Sci. Hungar.} {\bf 28} (1976), 77--80.

} 

\end{thebibliography}
 \end{document}